\documentclass[review]{elsarticle}
\usepackage{latexsym}
\usepackage{float}
\usepackage{pstricks}
\usepackage{amssymb}
\usepackage{amsmath}
\usepackage{amsfonts}
\usepackage{mathrsfs}

\usepackage{amsthm,array}
\usepackage{array}

\usepackage{graphicx}

\def\Ker{\text{\rm Ker\,}}

\def\dim{\text{\rm dim\,}}

\def\s1deg{S^1\text{\rm -Deg\,}}

\newcommand\bbZ{\ensuremath{\mathbb Z}}

\newrgbcolor{violet}{.6 .1 .8}
\newrgbcolor{lightyellow}{1 1 .8}
\newrgbcolor{lightblue}{.80 1 1}
\newrgbcolor{mygreen}{0 .66 .05}
\definecolor{mygreen}{rgb}{0,.66,.05}
\definecolor{lightyellow}{rgb}{1,1,.80}
\newrgbcolor{orange}{1 .6 0}
\newrgbcolor{GreenYellow}{.85 1 .31}
\newrgbcolor{Yellow}{1  1  0}
\newrgbcolor{Goldenrod}{1  .90  .16}
\newrgbcolor{Dandelion}{1  .71  .16}
\newrgbcolor{Apricot}{1  .68  .48}
\newrgbcolor{Peach}{1  .50  .30}
\newrgbcolor{Melon}{1  .54  .50}
\newrgbcolor{YellowOrange}{1  .58  0}
\newrgbcolor{Orange}{1  .39  .13}
\newrgbcolor{BurntOrange}{1  .49  0}
\newrgbcolor{Bittersweet}{1.  .4300  .24}
\newrgbcolor{RedOrange}{1  .23  .13}
\newrgbcolor{Mahogany}{1.  .4475  .4345}
\newrgbcolor{Maroon}{1.  .4084  .5376}
\newrgbcolor{BrickRed}{1.  .3592  .3232}
\newrgbcolor{Red}{1  0  0}
\newrgbcolor{OrangeRed}{1  0  .50}
\newrgbcolor{RubineRed}{1  0  .87}
\newrgbcolor{WildStrawberry}{1  .04  .61}
\newrgbcolor{Salmon}{1  .47  .62}
\newrgbcolor{CarnationPink}{1  .37  1}
\newrgbcolor{Magenta}{1  0  1}
\newrgbcolor{VioletRed}{1  .19  1}
\newrgbcolor{Rhodamine}{1  .18  1}
\newrgbcolor{Mulberry}{.6668  .1180  1.}
\newrgbcolor{RedViolet}{.9538  .4060  1.}
\newrgbcolor{Fuchsia}{.5676  .1628  1.}
\newrgbcolor{Lavender}{1  .52  1}
\newrgbcolor{Thistle}{.88  .41  1}
\newrgbcolor{Orchid}{.68  .36  1}
\newrgbcolor{DarkOrchid}{.60  .20  .80}
\newrgbcolor{Purple}{.55  .14  1}
\newrgbcolor{Plum}{.50  0  1}
\newrgbcolor{Violet}{.98 .15 .95}
\newrgbcolor{RoyalPurple}{.25  .10  1}
\newrgbcolor{BlueViolet}{.84  .38  .98}
\newrgbcolor{Periwinkle}{.43  .45  1}
\newrgbcolor{CadetBlue}{.38  .43  .77}
\newrgbcolor{CornflowerBlue}{.35  .87  1}
\newrgbcolor{MidnightBlue}{.4414  .9259  1.}
\newrgbcolor{NavyBlue}{.06  .46  1}
\newrgbcolor{RoyalBlue}{0  .50  1}
\newrgbcolor{Blue}{0  0  1}
\newrgbcolor{Cerulean}{.06  .89  1}
\newrgbcolor{Cyan}{0  1  1}
\newrgbcolor{ProcessBlue}{.04  1  1}
\newrgbcolor{SkyBlue}{.38  1  .88}
\newrgbcolor{Turquoise}{.15  1  .80}
\newrgbcolor{TealBlue}{.1572  1.  .6668}
\newrgbcolor{Aquamarine}{.18  1  .70}
\newrgbcolor{BlueGreen}{.15  1  .67}
\newrgbcolor{Emerald}{0  1  .50}
\newrgbcolor{JungleGreen}{.01  1  .48}
\newrgbcolor{SeaGreen}{.31  1  .50}
\newrgbcolor{Green}{0  1  0}
\newrgbcolor{ForestGreen}{.1992  1.  .2256}
\newrgbcolor{PineGreen}{.3100  1.  .5575}
\newrgbcolor{LimeGreen}{.50  1  0}
\newrgbcolor{YellowGreen}{.56  1  .26}
\newrgbcolor{SpringGreen}{.74  1  .24}
\newrgbcolor{OliveGreen}{.6160  1.  .4300}
\newrgbcolor{RawSienna}{.53  .28  .16}
\newrgbcolor{Sepia}{1.  .7510  .70}
\newrgbcolor{Brown}{.41  .25  .18}
\newrgbcolor{TAN}{.86  .58  .44}
\newrgbcolor{Gray}{1.  1.  1.}
\newrgbcolor{Black}{1  1  1}
\newrgbcolor{White}{1  1  1}

\newtheorem{theorem}{Theorem}[section]

\newtheorem{remark}[theorem]{Remark}

\newtheorem{remark-definition}[theorem]{Remark and Definition}

\usepackage{lineno,hyperref}
\modulolinenumbers[5]

\journal{Journal of differential equations}

\begin{document}

\begin{frontmatter}

\title{Non-invasive stabilization of periodic orbits in $O_4$-symmetrically coupled Van der Pol oscillators}


\author[utd]{Z. Balanov}
\ead{balanov@utdallas.edu}
\author[utd]{E. Hooton\corref{cor1}}
        \ead{exh121730@utdallas.edu}
\author[utd,gu]{W. Krawcewicz}  
\ead{wieslaw@utdallas.edu}
\author[utd]{D. Rachinskii}
\ead{dxr124030@utdallas.edu}
\address[utd]{Department of Mathematical Sciences,
        University of Texas at Dallas, Richardson, Texas, 75080 USA}
\address[gu]{College of Mathematics and Information Sciences, Guangzhou University,
Guangzhou, 510006
China}        
 \cortext[cor1]{Corresponding author}      
     


\begin{abstract}
Pyragas time delayed feedback control has proven itself as an effective tool to non-invasively stabilize periodic solutions. In a number of publications, this method was adapted to equivariant settings and applied to stabilize branches of small periodic solutions in systems of symmetrically coupled Landau oscillators near a Hopf bifurcation point. 
The form of the control ensures the non-invasiveness property, hence reducing the problem to  
finding a set of the gain matrices, which would guarantee the stabilization.
In this paper, we apply this method to a system of Van der Pol oscillators coupled in a cube-like configuration leading to $O_4$-equivariance. 
We discuss group theoretic restrictions which help to shape our choice of control.  
Furthermore, we explicitly describe the domains in the parameter space for which the periodic solutions are stable. 
\end{abstract}


\end{frontmatter}

\section{Introduction}
Stabilization of unstable periodic solutions is a classical control problem.
A control is called {\em non-invasive} if the controlled system has the same periodic solution as the uncontrolled system.
An elegant method of non-invasive control due to Pyragas \cite{pyragas1992continuous} is based on using a delayed phase variable.
This control strategy suggests to transform an uncontrolled ordinary differential system
\begin{equation}\label{1}
\dot{x}=F(x),\qquad x\in\mathbb{R}^N,
\end{equation}
into the delayed differential system
\begin{equation}\label{2}
\dot{x}=F(x) + {\mathcal K} (x(t)-x(t-\tau)),
\end{equation}
where ${\mathcal K}$ is the gain matrix. Obviously, if the delay $\tau$ equals the period $T_*$ of a periodic solution $x_*=x_*(t)$ to equation \eqref{1},
then $x_*$ is simultaneously a solution of the controlled equation \eqref{2}, since the control term 
${\mathcal K} (x(t)-x(t-\tau))$ vanishes on such solution. At the same time, Floquet multipliers of $x_*$ are different for the delayed and non-delayed equations, which may allow for stabilization
with a proper choice of the gain matrix $\mathcal{K}$. Note that typically the period $T_*$ of $x_*$ is not known {\em a priori}.
However, a stable periodic solution $x$ to \eqref{2} can usually be obtained for a range of delays $\tau$ sufficiently close to $T_*$.
Further tuning of the delay until the period $T$ of $x$ coincides with the delay $\tau$ can be used to achieve the non-invasive control.

A modification of the Pyragas control method that adapts it to symmetric (equivariant) setting has been developed in \cite{fiedler_Z_2,mamanya1,schneider_Z_3,schneider_elim}.
Periodic solutions of symmetric systems come in orbits (generated by the action of the symmetry group $G$ of the system)
and can be classified according to their symmetric properties. In particular, every periodic solution is fixed by a specific 
subgroup $H$ of the full group $G \times S^1$ of spatio-temporal symmetries. The control strategy proposed in \cite{fiedler_Z_2,mamanya1,schneider_Z_3,schneider_elim} is {\em selective}
in the sense that it acts non-invasively on periodic solutions with a specified period and  symmetry group including a given element 
while deforming or eliminating other periodic solutions. 
As a simple example, a control 
\begin{equation}\label{oddcont}
{\mathcal K} (x(t)+x(t-\tau/2))
\end{equation}can be used for non-invasive stabilization of $\mathbb{Z}_2$-symmetric
anti-periodic solutions $x(t)=-x(t-T/2)$ of period $T=\tau$, but this control does not vanish on $\tau$-periodic solutions that are not anti-periodic.

Since, in general, stability analysis of periodic solutions to delay differential equations, based on the usage of Floquet theory, is not well explored, in \cite{Fiedler_odd} it was suggested that complete stability analysis can be performed in the case of periodic solutions born via Hopf bifurcation.
Following \cite{Fiedler_odd} stability analysis for systems of symmetrically coupled Landau oscillators (the Landau oscillator is equivalent to the normal form of the Hopf bifurcation) was carried out in 
 \cite{fiedler_Z_2,mamanya1,schneider_Z_3,schneider_elim}, and essentially exploits the  idea outlined below. To conclude stability of a  bifurcating branch of periodic solutions from the well-known exchange of stability results, it is enough to check the following three conditions:

\begin{itemize}
\item[(i)] At the bifurcation point, the equilibrium is neutrally stable with neutral dimension two (genericity); 

\item[(ii)] The purely imaginary eigenvalues of the equilibrium cross the imaginary axis transversally;

\item[(iii)] The branch bifurcates in the direction in which the equilibrium becomes unstable. In the case of coupled Landau oscillators, this is simple to check since the periodic solutions are explicitly given.
\end{itemize}

For large dimensional delayed systems, stability of the equilibrium can be difficult to verify.
However, any application of the above control strategy to a specific symmetric system relies on the choice of one or several gain matrices. 
Since there is no general recipe for constructing those matrices, a possible approach to simplify analysis is to select a class of matrices depending on a small number of parameters. In particular, one can attempt to use {\em diagonal} (or block-diagonal) gain matrices, which allows one to  factorize the characteristic quasipolynomial.

In equivariant settings, it is usually the case that the  genericity condition (i) is violated. In \cite{fiedler_Z_2,mamanya1} the control \eqref{oddcont} is generalized to 
\begin{equation}
\label{int_cont}{\mathcal K} (\mathcal T_gx(t-2\pi\theta \tau)-x(t)),
\end{equation}
where $\mathcal T_g$ is the matrix associated with the {\it single} spatial group element $g$ and $2\pi\theta\tau$ is a rational fraction of the period which ``compensates'' the action of $g$ on the selected periodic solution.
Control \eqref{int_cont} 
breaks the symmetry (in particular, it is not non-invasive on the whole orbit of the targeted UPO)  and makes (i) possible to achieve for the controlled system.
At the same time, as was highlighted in \cite{mamanya1}, for certain groups, (i) can never be achieved by  \eqref{int_cont}. On the other hand, \cite{Isabelle_thesis} suggested a  general class of selective non-invasive equivariant Pyragas controls by taking a linear combination of controls of form \eqref{int_cont}, where $(g,\theta)$ varies amongst several group elements. 

In this paper, as a case study, we consider Hopf bifurcations in a system of 8 coupled Van der Pol oscillators arranged in a cubic connection
with a relatively complex group of permutational symmetries, $G=\mathbb{Z}_2 \times O_4$.
This system 
possesses one stable and 55 unstable branches of periodic solutions, which emanate from the zero equilibrium at 4 bifurcation points as a bifurcation parameter $\alpha$ is varied.
The branches can be classified into 12 types of spatio-temporal symmetries, which have been described in \cite{Coupled_Hys} using the equivariant topological degree method presented in  \cite{Green_book}.
We adapt one class of controls presented in \cite{Isabelle_thesis} with the objective to stabilize small periodic solutions from each 
branch using a selective control with the corresponding symmetry. 
We consider {\it linear combinations} of controls \eqref{int_cont} where we choose the values of $(g_k,\theta_k)$ from the symmetry group of the targeted unstable periodic solution in such a way that $\theta_k$ is constant. Also, we choose each ${\mathcal K}_k$ to be the same {\it real scalar matrix} with one scalar tuning parameter---the control strength $b$; another parameter is the coupling strength $a$ in the uncontrolled system. It turns out that these controls are sufficient for stabilizing unstable branches of all symmetry types except for one. Moreover, we obtain explicit expressions for stability domains in the $(a,b)$-plane for each stabilizable branch.
In Remark \ref{D_3}, we discuss a group-theoretic obstruction to this method and how this affects the branches which the chosen control fails to stabilize (cf. \cite{mamanya1}).

Unlike systems of Landau oscillators, the system of Van der Pol oscillators does not yield an explicit expression for periodic solutions in the 
form of relative equilibria. However, this does not create extra difficulties, since the proofs are based on asymptotic analysis. The proofs follow the general scheme from \cite{Fiedler_odd}.

The paper is organized as follows. In the next section, we describe symmetries of branches of periodic solutions for the system of interest and establish that all the four Hopf bifurcations giving rise to these branches are supercritical. Main results on stabilization of unstable branches by selective equivariant delayed control are presented in Section \ref{controlled}. Sections \ref{proofs} and \ref{conclusions} contain proofs and conclusions. The symbols representing spatio-temporal symmetry groups are explained in the Appendix.

\section{Uncontrolled system}\label{uncontrolled}

In this paper, we consider the system of coupled Van der Pol oscillators
\begin{equation} \label{vdpsystem}
\ddot x =  (\alpha -x^2) \dot x - x + \frac{a}{2}\mathcal B\dot x,
\end{equation}
where $x\in W := \mathbb{R}^8$; \ $\alpha$ is the bifurcation parameter, and
the interaction matrix has the form
\footnote{The system from  \cite{Coupled_Hys} describing an electrical circuit of coupled oscillators can be reduced to \eqref{vdpsystem} by standard rescaling. }
$$
{\mathcal B}=\left(
\begin{array}{cccccccc}
-3 & 1 & 0 & 1 & 1 & 0 & 0 & 0\\
1 & -3 &1 & 0  & 0 & 1 & 0 & 0\\
0 & 1 & -3 & 1 & 0 & 0 & 1 & 0\\
1 & 0 &1 & -3 & 0 & 0 & 0 & 1\\
1 & 0 & 0 & 0 & -3 & 1 & 0 & 1\\
0 & 1 & 0 & 0 & 1 & -3 & 1 & 0\\
0 &  0 &1 & 0 & 0 & 1 & -3 & 1\\
0 &  0 & 0 & 1 & 1 & 0 & 1 & -3
\end{array}
\right).
$$
The parameter $a$ measures the coupling strength.
In what follows, $V := W \oplus W$ stands for the phase space; also, the notation $x^3=x\cdot x\cdot x$ is used for componentwise multiplication, $(x\cdot y)_j=x_j y_j$, $j=1,\ldots,8$. 
For future reference, we denote the right hand side of \eqref{vdpsystem} by $f(\alpha,a,x,\dot x)$, which will allow us to use the notation
\begin{equation}\label{eq_abstract_eq}
\ddot{x}=f(\alpha,a,x,\dot x).
\end{equation}

In \cite{Coupled_Hys}, system \eqref{vdpsystem} was treated as an $S_4$-equivariant system, where $S_4 < S_8$ 
is the group of permutational symmetries of the cube preserving the 
orientation. If we include orientation reversing symmetries of the cube, this increases to $O_4 =  S_4 \times \mathbb Z_2$.  Noticing also that the right hand side of \eqref{vdpsystem} is an odd function (i.e. it is equivariant with respect to $\mathbb Z_2$ acting antipodally), in this paper we consider system \eqref{vdpsystem} with the full symmetry group $\mathbb Z_2 \times O_4 $. Each element $(r,g) \in   \mathbb Z_2 \times O_4$ is composed of $r=\pm 1$ and a permutation $g$ of $8$ symbols. We will denote by
$\mathscr T_g: W\to W$
the permutation matrix of $g$.

The spatio-temporal symmetries of a periodic function $x(t)$ are described by a subgroup $H < \mathbb Z_2 \times O_4$ and a homomorphism $\varphi: H \to S^1 \simeq \mathbb R /\mathbb Z$. This information is encoded in the graph of the homomorphism $\varphi$ which we will denote by $H^\varphi$. Put plainly, if $x(t)$ is a periodic function with period $T$ and symmetry group $H^\varphi$, then for each $(r,h) \in H$,
\begin{equation}\label{symm}
r \mathscr T_h x(t-\varphi(r, h)T) =  x(t).
\end{equation}

As it was shown in \cite{Coupled_Hys}, system \eqref{vdpsystem} undergoes 4 equivariant Hopf bifurcations giving rise to at least  56 branches of periodic solutions exhibiting different symmetry properties. Combining this result with the additional symmetry mentioned above allows us to describe the full symmetries of each branch (see Table 1 and Appendix for an explicit description of the groups listed in the second column).

\medskip

\begin{table}\label{table1}
\caption{Branches of solutions at Hopf bifurcation points}\label{11}
\begin{center}
\begin{tabular}{|c|c|c|}
\hline
Bifurcation~point & The group $H^\varphi$ of spatio-temporal & Total number of branches\\
& symmetries of the branch & at the bifurcation point $\alpha$\\
\hline
$\alpha = 0 $  & ${}^+(S_4)$ & 1 \\
$\alpha = a$ & $({}^-D^z_4), \; ({^-}D^z_3), \; ({}^-D^d_2), \:  ({}^-Z^c_4), \; ({}^-Z^t_3) $ & 27\\
$\alpha = 2a$& $(^+D^d_4), \; (^+D_3), \; (^+D^d_2), \; (^+Z^c_4), \; (^+Z^t_3)$  &  27\\
$\alpha = 3a$  & $(^-S^-_4 )$ & 1\\
\hline
\end{tabular}
\end{center}
\end{table}

To illustrate the meaning of these symmetries, let us take as an example the group
\begin{align*}
{}^-\mathbb Z^t_3 = & \{(1,(),0),(1,(245)(386),1/3),(1,(254)(368),2/3),\\ &(-1,(17)(28)(35)(46),0),(-1,(17)(265843),1/3),\\&(-1,(17)(234856),2/3),(-1,(),1/2),(-1,(245)(386),5/6),\\&(-1,(254)(368),1/6),(1,(17)(28)(35)(46),1/2),\\&(1,(17)(265843),5/6),(1,(17)(234856),1/6)\}.
\end{align*}
Suppose that $x(t)$ is a $T$-periodic function admitting the spatio-temporal symmetry ${}^-\mathbb Z^t_3$. Then, 
the components of $x$ respect certain relations. 
For example, 
for the element $(r,h,\varphi(r,h))=(-1,(17)(265843),1/3)\in {}^-\mathbb Z^t_3$, we have
$$
\mathscr T_h =  \left(
\begin{array}{cccccccc}
0& 0& 0 & 0 & 0 & 0 & 1 & 0 \\
0& 0& 1 & 0 & 0 & 0 & 0 & 0 \\
0& 0& 0 & 1 & 0 & 0 & 0 & 0 \\
0& 0& 0 & 0 & 0 & 0 & 0 & 1 \\
0& 0& 0 & 0 & 0 & 1 & 0 & 0 \\
0& 1& 0 & 0 & 0 & 0 & 0 & 0 \\
1& 0& 0 & 0 & 0 & 0 & 0 & 0 \\
0& 0& 0 & 0 & 1 & 0 & 0 & 0
\end{array}
\right)
$$
and, according to
\eqref{symm}, 
\begin{align*}
x_1(t)=-x_7(t-T/3)=x_1(t-2T/3)=&-x_7(t) =x_1(t-T/3)=-x_7(t-2T/3),\\
x_2(t)=-x_6(t-T/3)=x_5(t-2T/3)=&-x_8(t)=x_4(t-T/3)=-x_3(t-2T/3).
\end{align*}

The following statement plays an important role for the control problem.

\begin{theorem}\label{thm1}
All branches described in Table 1 are born via {\em supercritical} Hopf bifurcations.
\end{theorem}

The proof can be obtained by combining a standard asymptotic argument with $H$-fixed point reduction. We sketch the proof for convenience of the reader.

\begin{proof}
Notice that due to equivariance,  the space of $H$-fixed points 
$$
V^H:= \{x \in V \; : \; hx = x\;\; \forall  h \in H\}
$$ is a flow invariant subspace of the phase space for any $H < \mathbb Z_2 \times O_4$. If $x^*$ is a periodic solution with symmetry group $H^\varphi$, then $x^*(t) \in V^{{}_0H^\varphi}$ for all $t$ (cf. \cite{Green_book}), where  
\begin{equation}\label{def_0_H}
{}_0H^\varphi= \Ker \varphi.
\end{equation}
For each $H^\varphi$ appearing in Table 1, system \eqref{vdpsystem} restricted to $V^{{}_0H^\varphi}$  undergoes a (non-equivariant) Hopf bifurcation whose sub/supercriticality
coincides with that of the original system.   In what follows, we distinguish between the generic and non-generic (non-equivariant) Hopf bifurcations in the restricted systems.  
Our analysis splits into 3 cases when the Hopf bifurcation is generic and one case related to the non-generic setting.  

{\bf Case 1:} $H^\varphi = {}^+S_4$, $^-D_4^z$, $^-D_2^d$, $^+D_4^d$, $^+D_2^d$ or $^-S_4^-$. In this case, $V^{{}_0H^\varphi} = \mathbb R \oplus \mathbb R$ and system \eqref{vdpsystem} restricted to $V^{{}_0H^\varphi}$ is equivalent to the equation of a single Van der Pol oscillator, where the parameter is shifted by an integer multiple of $a$ depending on the branch.

\smallskip
{\bf Case 2:} $H^\varphi = {}^-\mathbb Z_4^c$ or $^+\mathbb Z_4^c$. In this case, $V^{{}_0H^\varphi}=\mathbb R^2 \oplus \mathbb R^2$ and system \eqref{vdpsystem} restricted to $V^{{}_0H^\varphi}$ is equivalent to the system  of two uncoupled Van der Pol oscillators. This system has a continuum of periodic solutions depending on the phase between the two oscillators. Although they all correspond to solutions of the original system, only the solutions for which the first oscillator is one quarter of the period out of phase with the second correspond to the solutions with the prescribed symmetry.

\smallskip
{\bf Case 3:}  $H^\varphi = {}^+D_3$ or $^-D_3^z$.  In this case, $V^{{}_0H^\varphi}=\mathbb R^2 \oplus \mathbb R^2$ and system \eqref{vdpsystem} restricted to $V^{{}_0H^\varphi}$ is equivalent to the system of two asymmetrically coupled Van der Pol oscillators given by
\begin{align*}
{}^-D_3: \quad \begin{cases}
\ddot x_1-\alpha \dot x_1+\dot x_1x^2_1+x_1=\frac{3a}{2}(\dot x_2-\dot x_1),\\
\ddot x_2-\alpha \dot x_2+\dot x_2x^2_2+x_2=\frac{a}{2}(\dot x_1-5\dot x_2);
\end{cases}
\\
{}^+D_3^z: \quad \begin{cases}
\ddot x_1-\alpha \dot x_1+\dot x_1x^2_1+x_1=\frac{3a}{2}(\dot x_2-\dot x_1),\\
\ddot x_2-\alpha \dot x_2+\dot x_2x^2_2+x_2=\frac{a}{2}(\dot x_1-\dot x_2).
\end{cases}
\end{align*}

{\bf Case 4:} $H^\varphi = {}^-\bbZ_3^t$ or $^+\bbZ_3^t$. In this case, $V^{{}_0H^\varphi}=\mathbb R^4 \oplus \mathbb R^4$ and system \eqref{vdpsystem} restricted to $V^{{}_0H^\varphi}$ undergoes a non-generic Hopf bifurcation. We will consider the following families of non-symmetric delayed differential equations with an additional parameter   $T$:
\begin{align}\label{eq:reduced_Z_3}
{}^-\mathbb Z_3^t: \quad \begin{cases}
\ddot x_1-\alpha \dot x_1+\dot x_1x^2_1+x_1=\frac{a}{2}\left(\dot x_2+\dot x_2(t-\frac{T}{3})+\dot x_2(t-\frac{2T}{3})-3\dot x_1\right),\\
\ddot x_2-\alpha \dot x_2+\dot x_2x^2_2+x_2=\frac{a}{2}\left(\dot x_1-\dot x_2(t-\frac{T}{3})-\dot x_2(t-\frac{2T}{3})-3\dot x_2\right);
\end{cases}\\ \label{eq_reduced_+Z_3}
{}^+\mathbb Z_3^t: \quad \begin{cases}
\ddot x_1-\alpha \dot x_1+\dot x_1x^2_1+x_1=\frac{a}{2}\left(\dot x_2+\dot x_2(t-\frac{T}{3})+\dot x_2(t-\frac{2T}{3})-3\dot x_1\right),\\
\ddot x_2-\alpha \dot x_2+\dot x_2x^2_2+x_2=\frac{a}{2}\left(\dot x_1+\dot x_2(t-\frac{T}{3})+\dot x_2(t-\frac{2T}{3})-3\dot x_2\right).
\end{cases}
\end{align}
Clearly, $T$-periodic solutions to the original system with the spatio-temporal symmetry $ {}^-\bbZ_3^t$  (resp. $^+\bbZ_3^t$) are in one-to-one correspondence
with $T$-periodic solutions to \eqref{eq:reduced_Z_3}  (resp. \eqref{eq_reduced_+Z_3}).

To establish  supercriticality, in  Cases $1$ and $2$ we recall that the branch of periodic solutions of a Van der Pol equation is supercritical, while  for Cases $3$ and $4$ one can apply the standard techniques of asymptotic analysis. We will just give a detailed explanation for \eqref{eq:reduced_Z_3} since the other cases are analogous.

{\bf Step 1:} By rescaling time $y(\beta t)=x(t)$, where $\beta = T/2\pi$, one obtains:
\begin{align*}
\beta^2\ddot y_1-\alpha \beta \dot y_1+\beta \dot y_1y^2_1+y_1=\frac{a\beta}{2}\left(\dot y_2+\dot y_2\left(t-\frac{2\pi}{3}\right)+\dot y_2\left(t-\frac{4\pi}{3}\right)-3\dot y_1\right),\\
\beta^2 \ddot y_2-\alpha \beta \dot y_2+\beta \dot y_2y^2_2+y_2=\frac{a\beta}{2}\left(\dot y_1-\dot y_2\left(t-\frac{2\pi}{3}\right)-\dot y_2\left(t-\frac{4\pi}{3}\right)-3\dot y_2\right).
\end{align*}

{\bf Step 2:} We will take $r$ to be a small parameter and expand the parameters $\alpha$ and $\beta$ near the values $\alpha = a$ and $\beta = 1$  as follows:
\begin{align*}
\alpha = a + \hat{\alpha}r^2 +o(r^2), \qquad
\beta= 1 + \hat{\beta}r^2 +o(r^2).
\end{align*}
The standard results about asymptotics of branches born at a Hopf point legitimize the absence of linear terms. We will now expand 
$$ 
y_2 =r\cos t + r^3 \psi_2(t) +o(r^3), 
$$
where $\psi_2$ is orthogonal to $\sin t$ and $\cos t$ in $\mathbb{L}_2[0,2\pi]$. 
Plugging this expression into the first equation shows that $y_1$ has only harmonics of order divisible by $3$ and its expansion starts with $r^3$. This allows us to expand
$$
y_1 =r^3 \psi_1(t)+o(r^3),
$$
where $\psi_1(t)$ is orthogonal to $\sin t$ and $\cos t$.

{\bf Step 3:} Projecting terms of order $r^3$ in the second equation onto the first Fourier mode gives the equation
$$
-2\hat\beta\cos t +(a\hat{\beta}+\hat{\alpha})\sin t -\frac{1}{4}\sin t = a\hat{\beta}\sin t.
$$
From this it can be seen that $\hat\alpha= 1/4 > 0$, so the branch must be supercritical.
\end{proof}

\medskip
\begin{remark}\label{rem:supercriticality}
{\rm Theorem \ref{thm1} allows us to reduce the analysis of stability of periodic solutions to studying characteristic equations related to the zero equilibrium, from which the periodic solutions bifurcate.}
\end{remark}

\section{Main Results}\label{controlled}
For the symmetry group $H^\varphi$, recall that ${}_0H^\varphi= \ker \varphi$ (cf. \eqref{def_0_H}). We will denote by $t_0(H^\varphi)$ 
the smallest $t \in (0,1)$ such that $t = \varphi(r,h)$ for some $(r,h) \in H$. Finally, define a set of spatial symmetries by
$${}_1H^\varphi= \varphi^{-1}(t_0(H^\varphi))$$
and by $|H|$ the cardinality of $H$.

\begin{theorem}\label{thm2}
Suppose $x^*_{\alpha}$ is a branch of periodic solutions to \eqref{eq_abstract_eq}  with symmetry $K^\varphi = {}^-D_4^z$, $^-D_2^d$, ${}^-D_3^z$, $^+D_4^d$, $^+D_2^d$ or $^-S_4^-$ which bifurcates from the zero solution $x=0$ at $\alpha_o = ka$ (where $k = 1,2,3$ is given in Table \ref{11}). Then, for every $b > ka$ there exists an $\alpha^ *=\alpha^*(a,b) > \alpha_o$ such that $x^*_{\alpha}$ is an asymptotically stable solution of
\begin{equation}\label{eq_cont_eq1}
\ddot{x}= f(\alpha,a,x,\dot x) +  b\left(-\dot x(t)+\frac{1}{|{}_0H^\varphi |}\sum_{(r,h)\in {}_0H^\varphi}r\mathscr T_h \dot x(t)\right)
\end{equation}
for every $\alpha \in (\alpha_o,\alpha^*)$.
\end{theorem}

\begin{theorem}\label{thm3}
Suppose $x^*_{\alpha}$ is a branch of $T_{\alpha}$-periodic solutions to \eqref{eq_abstract_eq}  with symmetry $H^\varphi = {}^-\mathbb Z_4^c$, $^-\mathbb Z_3^t$, $ {}^+\mathbb Z_4^c$ or $^+\mathbb Z_3^t$  which bifurcates from $x=0$ at $\alpha_o = ka$ (where $k = 1,2$ is given in Table \ref{11}). Then, there exists a domain $\mathcal D\in\mathbb R^2_+$ such that for every point $(a,b)\in \mathcal D$ there exists an $\alpha^* = \alpha^*(a,b) > \alpha_o$ such that  $x^*_{\alpha}$ is an asymptotically stable solution of
\begin{equation}\label{eq_cont_eq2}
\ddot{x}= f(\alpha,a,x,\dot x) +  b\left(-\dot x(t)+\frac{1}{|{}_1H^\varphi |}\sum_{(r,h)\in {}_1H^\varphi}r\mathscr T_h \dot x(t-\tau_\alpha)\right).
\end{equation}
for every $\alpha \in (\alpha_o ,\alpha^*)$ with $\tau_\alpha=t_0(H^\varphi) T_{\alpha}$. Furthermore, for each $H^\varphi$, the domain $\mathcal D$  is  explicitly  described in Table \ref{22}.
\end{theorem}

\begin{table}\label{table2}
\caption{Domains of stability}\label{22}
\begin{center}
\begin{tabular}{|c|c|}
\hline
Symmetry of the branch & Domain $\mathcal D$ of parameters for which the branch is stable\\
\hline
$({}^-\mathbb Z^c_4)$ & $0<a<b$ \\
$({}^-\mathbb Z^t_3)$ & $0<a<b$\\
$({}^+\mathbb Z^c_4)$ & $0<2a<b$ \\
$({}^+\mathbb Z^t_3)$ & $0<a<\psi(b)$, where $\psi$ is described in \\ & Remark \ref{YTM} and illustrated in Figure \ref{siska} \\
\hline
\end{tabular}
\end{center}
\end{table}

\begin{figure}
\begin{center}
\includegraphics*[width=0.5\columnwidth]{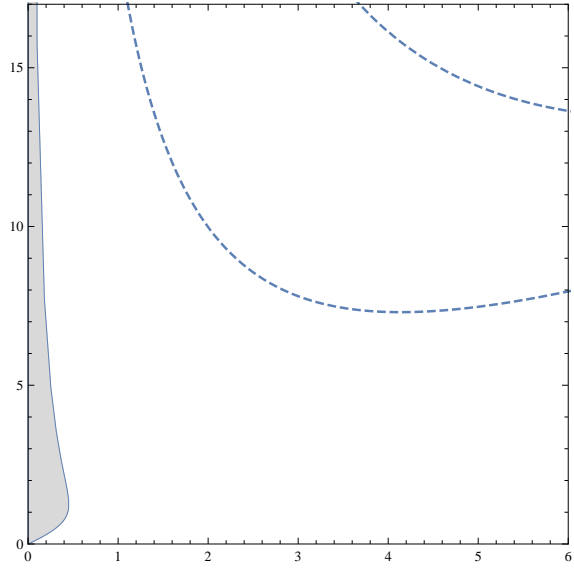} 
\end{center}
\caption{Domain ${\mathcal D}$ of stability of the branch with symmetry ${}^+\mathbb Z^t_3$ on the $(a,b)$-plane (shaded).
}
\label{siska}
\end{figure}

\begin{remark}\label{YTM}
{\rm Consider the
curve 
\begin{equation}\label{xxx}
(a,b)=(\gamma_1(s),\gamma_2(s))=\left(\frac{(s^2 - 1) (1 + \cos(\frac{s\pi}{3}))}{2s\sin(\frac{s\pi}{3})},\frac{s^2 - 1}{ s \sin(\frac{s\pi}{3})} \right), \qquad s\in [1,3),
\end{equation}
which bounds the shaded domain ${\mathcal D}$ in Figure \ref{siska}.
By direct computation, it is easy to see that $\gamma_2(s)$ is monotonic on the interval $[1,3)$, and therefore invertible. The function $\psi$ appearing in Table \ref{22} is defined by $\psi:= \gamma_1\circ\gamma_2^{-1}$.}
\end{remark}

\begin{remark}\label{D_3}
{\rm Since our analysis of stability of the bifurcating branch (with symmetry $H^\varphi$) in the controlled system relies on the standard exchange of stability results, we require that $\pm i$ has multiplicity one at the bifurcation point. 
For an element $(r,h,\theta) \in H^\varphi$ we will denote by $V^c_{(r,h,\theta)}$ the set of points in the complexification of center space which is fixed by $(r,h,\theta)$, where $\theta$ acts on the complexification by multiplication by $e^{2\pi i \theta}$. 
It was observed in \cite{mamanya1} that for a control of the form $\mathcal K(r\mathcal T_hx(t-2\pi\theta T)-x(t))$, the above condition can be satisfied only if  $\dim_\mathbb C V^c_{(r,h,\theta)} = 2$.
For any subset $S:=\{(r_k,h_k,\theta_k)\}$ of $H^\varphi$, the equivalent requirement for a linear combination of these controls is that $$\dim_\mathbb C \bigcap_{(r_k,h_k,\theta_k)\in S} V^c_{(r_k,h_k,\theta_k)} = 2.$$
For the majority of the branches considered in this paper, although for a single group element $(r,h,\theta)$ this condition is not satisfied, it is satisfied if we consider a set $S$ of several group elements where $\theta_k\equiv \theta$ is the same for all $k$ (cf. \eqref{eq_cont_eq2}).
However, in the case of ${}^+D_3$ and $a=0$, we have  
$$\dim_\mathbb C \bigcap_{(r,h,\theta)\in{}^+D_3} V^c_{(r,h,\theta)} = 4.$$
For this reason, our control fails to stabilize the branch with symmetry ${}^+D_3$.
It is our conjecture (confirmed by numerical simulations) that this obstruction still exists for weak coupling. 
}
\end{remark}

\section{Proofs}\label{proofs}
The $\mathbb Z_2 \times O_4$-isotypical decomposition of $V = W \oplus W$ is given by
\begin{equation}
W = W_1 \oplus W_2 \oplus W_3 \oplus W_4,
\end{equation}
where $W_1$ (resp. $W_2$, $W_3$ and $W_4$) are mutually non-equivalent absolutely irreducible representations with $\dim W_1 = 1$ (resp.  $\dim W_2 = 3$, $\dim W_3 = 3$ and 
$\dim W_4 = 1$).
 Take a basis $e_1 \in W_1$ (resp. $e_2$, $e_3$, $e_4 \in W_2$, $e_5$, $e_6$, $e_7 \in W_3$, $e_8 \in W_4$) and call the basis $e_1, \ldots, e_8 \subset W$ an {\it isotypical basis} for $W$. Observe that in any isotypical basis 
the linearization of system \eqref{vdpsystem} at the origin is given by
\begin{equation}\label{eq_lin_iso}
\ddot x = A_0 \dot x - x
\end{equation}
with
\begin{equation}\label{a0}
A_0=\left(
\begin{array}{cccccccc}
\alpha& 0& 0 & 0 & 0 & 0 & 0 & 0 \\
 0& \alpha-a& 0 & 0 & 0 & 0 & 0 & 0 \\
 0& 0&\alpha-a & 0 & 0 & 0 & 0 & 0 \\
0& 0& 0 & \alpha-a & 0 & 0 & 0 & 0 \\
0& 0& 0 & 0 & \alpha-2a & 0 & 0 & 0 \\
0& 0& 0 & 0 & 0 & \alpha-2a & 0 & 0 \\
0& 0& 0 & 0 & 0 & 0 & \alpha-2a & 0 \\
0& 0& 0 & 0 & 0 & 0 & 0 & \alpha-3a
\end{array}
\right).
\end{equation}
Hereafter, we will assume that the linearized system is of the form \eqref{eq_lin_iso}.

\subsection{Proof of Theorem \ref{thm2}}
Since the treatment of each branch relevant to this theorem follows the same lines, we restrict ourselves to the case when $H^\varphi= {}^-D_2^d$
for which we have
\begin{align*}
{}^-_0D_2^d=& \big\{\big(1,()\big),\big(-1,(13)(24)(57)(68)\big),\big(1,(15)(28)(37)(46)\big),\big(-1,(17)(26)(35)(48)\big), \\& \big(-1,(17)(28)(35)(46)\big),\big(1,(15)(26)(37)(48)\big),\big(-1,(13)(57),\big),\big(1,(24)(68)\big)\big\}.
\end{align*}
Then, it follows that the control term (written in the {\it original basis}) is represented by
\begin{align}\label{eq_cont}
 b\left(-\dot x(t)+\frac{1}{|{}_0H^\varphi |}\sum_{h\in {}_0H^\varphi}\mathscr T_h \dot x(t)\right)= \frac{b}{4}\left(
\begin{array}{cccccccc}
-3 & 0& -1& 0& 1& 0& -1& 0\\
0& -4& 0& 0& 0& 0& 0 & 0\\
-1& 0& -3& 0& -1& 0& 1& 0\\
0& 0& 0& -4& 0& 0& 0& 0\\
1& 0& -1& 0& -3& 0& -1& 0\\
0& 0& 0& 0& 0& -4& 0& 0\\
-1& 0& 1& 0& -1& 0& -3&, 0\\
0& 0& 0& 0& 0& 0& 0& -4
\end{array}
\right)\dot x.
\end{align}
Notice that $x_\alpha$ bifurcates at the value $\alpha_o=a$ (see Table \ref{11}). Due to Theorem \ref{thm3}, to complete the proof, it is enough to show that if $b>a$, then the unstable dimension of the trivial equilibrium of system \eqref{eq_cont_eq1}  changes from zero  to two  as $\alpha$ increases and passes $\alpha_o$.  Combining \eqref{eq_lin_iso}   with \eqref{eq_cont} (written in the  {\it isotypical basis}) allows us to write the linearization of \eqref{eq_cont_eq1} as 
$$
\ddot x= (A_0-b B_0)\dot x-x
$$
with the matrix $A_0$ defined by \eqref{a0} and
$$
B_0=\left(
\begin{array}{cccccccc}
1& 0& 0 & 0 & 0 & 0 & 0 & 0 \\
 0& 0& 0 & 0 & 0 & 0 & 0 & 0 \\
 0& 0&1 & 0 & 0 & 0 & 0 & 0 \\
0& 0& 0 & 1 & 0 & 0 & 0 & 0 \\
0& 0& 0 & 0 & 1 & 0 & 0 & 0 \\
0& 0& 0 & 0 & 0 & 1 & 0 & 0 \\
0& 0& 0 & 0 & 0 & 0 & 1 & 0 \\
0& 0& 0 & 0 & 0 & 0 & 0 & 1
\end{array}
\right).
$$
Since  $b>a>0$, it is easy to see that for $\alpha$ close to $a$, all but one pair of eigenvalues have negative real part, while the real part of that pair increases as  $\alpha$ increases and passes $\alpha_o$. This completes the proof.

\subsection{Proof of Theorem \ref{thm3}} 
The proof of Theorem  \ref{thm3}  
requires that for each $H^{\varphi}$ one computes the characteristic equation of the linearization of system \eqref{eq_cont_eq2} 
at the origin. The results of these computations done in an isotypical basis are presented in Table \ref{33}. Since the treatment of each branch appearing in Table \ref{33} follows the same lines, we restrict ourselves to the case when $H^\varphi= {}^+\mathbb Z_3^t$. Similarly to the proof of Theorem \ref{thm2}, our goal is  
to show that if $(a,b) \in \mathcal D$, then the unstable dimension of the trivial equilibrium of system \eqref{eq_cont_eq2}  changes from zero  to two  as $\alpha$ increases 
and passes $\alpha_o$. This goal is achieved in two steps.

\medskip
{\bf Step 1.} At this stage we show that for $\alpha = \alpha_o$ and any $(a,b) \in \mathcal D$, the trivial equilibrium of system \eqref{eq_cont_eq2} has  a two-dimensional center manifold and no unstable manifold. To this end, taking characteristic equations from Table \ref{33} related to $H^\varphi = {}^+\mathbb Z^t_3$, and putting $\alpha = \alpha_o = 2a$ and $T_\alpha = 2 \pi$ yields the following equations (here $\eta = e^{{2\pi \over 6}i}$):
\begin{align} \label{eq_char1}
\lambda^2+(b-2a)\lambda +1 &= -b\lambda e^{-2\lambda\pi\over 6},\\ \nonumber
\lambda^2+(b-a)\lambda +1 &= 0,\\ \nonumber
\lambda^2+(b-a)\lambda +1 &= 0,\\ \nonumber
\lambda^2+(b-a)\lambda +1 &= 0,\\ \label{eq_char2}
\lambda^2+b\lambda +1 &= -b\lambda e^{-2\lambda\pi\over 6},\\ \label{eq_char3}
\lambda^2+b\lambda +1 &= b\eta \lambda e^{-2\lambda\pi\over 6},\\ \label{eq_char4}
\lambda^2+b\lambda +1 &= b \overline\eta \lambda e^{-2\lambda\pi\over 6},\\ \nonumber
\lambda^2+(a+b)\lambda +1 &= 0.
\end{align}
The spectrum of the zero equilibrium is the union of all the solutions $\lambda$ to these 8 equations.
By inspection, if $b>a>0$, then, except for \eqref{eq_char1}---\eqref {eq_char4}, the above equations do not admit roots with non-negative real parts meaning that the corresponding polynomials are stable.
Next, notice that for  $b=0$, equations \eqref{eq_char1}---\eqref{eq_char4} admit $\pm i$ as a root. Furthermore, for any $b$, $i$ remains a root of \eqref{eq_char3}, while $-i$ remains a root of \eqref{eq_char4}. Finally observe that, by implicit differentiation of equations \eqref{eq_char1}---\eqref{eq_char4}  with respect to $b$ at $a = b = 0$ and $\lambda = \pm i$, it is easy to see that for any sufficiently small $b >0$,  all other roots of  \eqref{eq_char1}---\eqref{eq_char4}  lie in the left half plane. 

\begin{table}\caption{Characteristic equations written in isotypical coordinates.}\label{33}
\begin{center}
\begin{tabular}{|c|c|c|c|c|c|}
\hline
&&${}^-\mathbb{Z}_4^c$&${}^-\mathbb{Z}_3^t$&${}^+\mathbb{Z}_4^c$&${}^+\mathbb{Z}_3^t$\\
\hline
$1$ &$\lambda^2+(b-\alpha)\lambda +1 =$ & $0$ & $0$ & $0$& $-b\lambda e^{\frac{-T\lambda}{6}}$\\
$2$ &$\lambda^2+(b+a-\alpha)\lambda +1= $ &$ib\lambda e^{\frac{-T\lambda}{4}}$&$-b\lambda e^{\frac{-T\lambda}{6}}$&$0$&$0$\\
$3$ & $\lambda^2+(b+a-\alpha)\lambda +1 = $ &$-ib\lambda e^{\frac{-T\lambda}{4}}$&$ b\lambda e^{\frac{2\pi i-T\lambda}{6}}$&$0$&$0$\\
$4$ &$\lambda^2+(b+a-\alpha)\lambda +1 = $ &$0$&$b\lambda  e^{\frac{-2\pi i-T\lambda}{6}}$&$0$&$0$\\
$5$ &$\lambda^2+(b+2a-\alpha)\lambda +1 = $ &$0$&$0$&$ib\lambda e^{\frac{-T\lambda}{4}}$& $-b\lambda e^{\frac{-T\lambda}{6}}$\\
$6$ &$\lambda^2+(b+2a-\alpha)\lambda +1 = $&$0$&$0$&$-ib\lambda e^{\frac{-T\lambda}{4}}$& $b\lambda e^{\frac{2\pi i-T\lambda}{6}}$\\
$7$ &$\lambda^2+(b+2a-\alpha)\lambda +1 = $&$0$&$0$&$0$&$b\lambda e^{\frac{-2\pi i-T\lambda}{6}}$\\
$8$ &$\lambda^2+(b+3a-\alpha)\lambda +1 = $&$0$&$-b\lambda e^{\frac{-T\lambda}{6}}$&$0$&$0$ \\
\hline
\end{tabular}
\end{center}
\end{table}

To show that for any $(a,b)\in \mathcal D$,  the roots of equations  \eqref{eq_char1}---\eqref{eq_char4} lie in the left half plane, we use a variant of  Zero Exclusion Principle. Since 
$\mathcal D$ contains the points $(0,b)$ for small $b > 0$, it is enough to show that as $(a,b)$ varies in $\mathcal D$, no roots of equations \eqref{eq_char1}---\eqref{eq_char4} 
can ever pass through the purely imaginary axis. To this end,  we plug $\lambda = is$ into each equation in turn. The points $(a,b)$ for which  \eqref{eq_char1} admits a purely imaginary root form the set of curves given by
\[
\gamma(s)=(a(s),b(s))=\left(\frac{(s^2 - 1) (1 + \cos(\frac{s\pi}{3}))}{2s\sin(\frac{s\pi}{3})},\frac{s^2 - 1}{ s \sin(\frac{s\pi}{3})} \right)
\]
(cf.~\eqref{xxx}).
Since for each $s$, $$\frac{a(s)}{b(s)}= \frac{1 + \cos(\frac{s\pi}{3})}{2}\le 1$$
with equality iff $s=6k+1$, for some integer $k$  the segment of the curve corresponding to $1<s<3$  lies above  the straight line $a=b$. Notice that if  $\gamma(s)=\gamma(t)$ for some $t \neq s$, then either $s=t+6k$ or $s= 6k-t$ for some integer $k$. Combining this observation with monotonicity of $s^2-1$ for $s>1$ and periodicity of the sinus function proves that $\gamma(s)$ does not have self-intersection points. For this reason, we see that $\mathcal D$ is bounded by $\{\gamma(s) : s\in [0,3)\}\cup \{ a=0\}$ (cf. Remark \ref{YTM} and Figure \ref{siska}).

By taking the absolute value of both sides of \eqref{eq_char2}, it follows that \eqref{eq_char2} never admits a purely imaginary root. On the other hand, while \eqref{eq_char3} admits $i$ as a root for all $(a,b)$, the same argument as for \eqref{eq_char2}  shows that $\lambda = i$ is the only purely imaginary root of \eqref{eq_char3}.
By differentiating \eqref{eq_char3} with respect to $\lambda$,  one concludes that for $b>0$, $\lambda = i$  is a {\it simple} root.
Replacing $i$ by $-i$, one can apply the same argument to \eqref{eq_char4}.

{\it To summarize Step 1: } We showed that, in the case of  ${}^-\mathbb Z_3^t$, at the bifurcation point $\alpha =2a$, if $a,b>0$, only one of the quasi-polynomials from Table \ref{33}  can admit roots with positive real part (namely equation $1$).
On the other hand, the boundary of $\mathcal D$ is defined by the values of $(a,b)$ for which equation $1$ admits purely imaginary roots.
In the cases of ${}^-\mathbb Z_4^c$, ${}^-\mathbb Z_3^t$ and ${}^+\mathbb Z_4^c$, at the corresponding bifurcation points, all the quasi-polynomials from Table  \ref{33} do not admit roots with positive real parts. 
This explains why the case of ${}^-\mathbb Z_3^t$ was taken as the  demonstrative example and why in Table \ref{22} it has a seemingly peculiar entry.

\medskip

{\bf Step 2.} It is now left to show that as $\alpha$ increases and $T(\alpha)$ varies, the purely imaginary root $i$ (resp. $-i$) of equation 6 (resp. 7) in Table \ref{33}
moves into the right half-plane. To this end, 
following the idea suggested in \cite{Fiedler_odd}, p. 326 (see also references therein), we will fix $a$ and $b$,  and treat $\alpha$ and $T$ as independent bifurcation parameters. 
Let us show that in the $(\alpha,T)$-plane a Hopf curve passes through the point $(2a, 2\pi)$ with a vertical tangent line. In fact, substituting $i \omega$ into Table \ref{33}, equation 6, one obtains
$$
1 - \omega^2 + (b + 2a - \alpha) i\omega = b\eta i\omega e^{-i\omega T \over 6}.
$$ 
The above equation implicitly defines $\alpha$ and $T$ as functions of $\omega$. Differentiating this equation with respect to $\omega$ and separating real and imaginary parts yields
\begin{equation}\label{eq:pots1}
-2 = {b \over 6} (2\pi + T^{\prime}(\omega)),\qquad
\alpha^{\prime}(\omega)  = 0
\end{equation}
as desired. On the other hand, fixing $T = 2\pi$ and differentiating equation 6 from Table \ref{33} with respect to $\alpha$ at  $\alpha = 2a$ and $\lambda = i$ yields: 
\begin{equation}\label{eq:pots2}
\lambda^{\prime} = \frac{3} {6 + \pi b} > 0.
\end{equation}
Combining \eqref{eq:pots1} and \eqref{eq:pots2} implies that for any function  $T=T(\alpha)$ with $T(2a)=2\pi$, one has that $\lambda^\prime$ evaluated at $\alpha = 2a$ and $\lambda= i$, is positive. The same argument can be used in the case of $-i$ as a root of Table \ref{33},  equation 7. Combining this with Theorem \ref{thm1} and the standard exchange of stability results completes the proof.

\section{Conclusions}\label{conclusions}
We have considered a system of symmetrically coupled Van der Pol oscillators with $O_4$-permu\-tational symmetry.
This system possesses multiple branches of unstable periodic solutions with different symmetry properties.
Using an equivariant Pyragas type delayed control introduced in \cite{fiedler_Z_2,mamanya1,schneider_Z_3} we proposed a specific form of the gain matrices, which ensures the non-invasive 
stabilization of periodic solutions near a Hopf bifurcation point for the branches of each symmetry type with one exception.
We found explicitly stability domains of the controlled system in the parameter space.
The failure of the control for branches with one specific type of symmetry can be associated 
with group theoretic restrictions considered in \cite{mamanya1}.

\section{Appendix}
In this Appendix, we will explain the symbols used in the main body of the text to denote spatio-temporal symmetry groups.  For any $H<S_4\times S^1$, we will define ${}^-H <\mathbb Z_2\times O_4\times S^1$ and ${}^+H < \mathbb Z_2\times O_4\times S^1$
by 
\begin{align*}
{}^+H :=H \times(\mathbb Z_2\times O_1)^{o},\qquad
{}^-H :=H \times(\mathbb Z_2\times O_1)^{oz},
\end{align*}
where 
$$(\mathbb Z_2\times O_1  )^{o} := \{\big(1,( ),0\big), \big(1,(17)(28)(35)(46),0\big), (-1,( ),1/2),  (-1,(17)(28)(35)(46),1/2,)\}$$
$$(\mathbb Z_2 \times O_1  )^{oz} := \{(1,( ),0), (-1,(17)(28)(35)(46),0), (-1,( ),1/2),  (1,(17)(28)(35)(46),1/2)\}.$$
All spatio-temporal symmetry groups which we deal with in this paper  appear as either $^+H$ or $^-H$, where $H$ is among the following groups:
\begin{align*}
S_4:=&\{((),0),((15)(28)(37)(46),0),((17)(26)(35)(48),0),((12)(35)(46)(78),0),((17)(28)(34)(56),0),\\&((14)(28)(35)(67)0),((17)(23)(46)(58),0),((13)(24)(57)(68),0),((18)(27)(36)(45),0),\\&((16)(25)(38)(47),0),((254)(368),0),((245)(386),0),((163)(457),0),((136)(475),0),((168)(274),0),\\&((186)(247),0),((138)(275),0),((183)(257),0), ((1234)(5678),0),((1432)(5876),0),\\&((1265)(3874),0),((1562)(3478),0),((1485)(2376),0),((1584)(2678),0)\}\\
D_4^z:=&\{((),0),((1234)(5678),0),((13)(24)(57)(68),0),((1432)(5876),0),\\&((17)(26)(35)(48),1/2),((18)(27)(36)(45),1/2),((15)(28)(37)(46),1/2),((16)(25)(38)(47),1/2)\}\\
D_3^z:=&\{((),0),((254)(368),0),((245)(386),0),((17)(26)(35)(48),1/2),\\&((17)(28)(34)(56),1/2),((17)(23)(46)(58),1/2)\}\\
D_2^d:=&\{((),0),((17)(26)(35)(48),1),((13)(24)(57)(68),1/2),((15)(28)(37)(46),1/2)\}\\
\mathbb Z^c_4:=&\{((),0),((1234)(5678),1/4),((13)(24)(57)(68),1/2),((1432)(5876),3/4)\}\\
\mathbb{Z}^t_3:=&\{((),0),((254)(368),1/3),((245)(386),2/3)\}\\
D_4^d:=&\{((),0),((1234)(5678),1/2),((13)(24)(57)(68),0),((1432)(5876),1/2),\\&((17)(26)(35)(48),0),((18)(27)(36)(45),1/2),((15)(28)(37)(46),0),((16)(25)(38)(47),1/2)\}\\
D_3:=&\{((),0),((254)(368),0),((245)(386),(17)(26)(35)(48),0),((17)(28)(34)(56),0),((17)(23)(46)(58),0)\}\\
S_4^-:=&\{((),0),((15)(28)(37)(46),1/2),((17)(26)(35)(48),1/2),((12)(35)(46)(78),1/2),((17)(28)(34)(56),1/2),\\&((14)(28)(35)(67)1/2),((17)(23)(46)(58),1/2), ((13)(24)(57)(68),0),((18)(27)(36)(45),0),\\&((16)(25)(38)(47),0),((254)(368),0),((245)(386),0),((163)(457),0),((136)(475),0),((168)(274),0),\\&((186)(247),0),((138)(275),0),((183)(257),0), ((1234)(5678),1/2),((1432)(5876),1/2),\\&((1265)(3874),1/2),((1562)(3478),1/2),((1485)(2376),1/2),((1584)(2678),1/2)\}
\end{align*}

\section*{Acknowledgements} The authors acknowledge the support from National Science Foundation through grant DMS-1413223. The first author is grateful for the support of the Gelbart Research Institute for mathematical sciences at Bar Ilan University. The third author was also supported by National Natural Science Foundation of China (no. 11301102).

\section*{References}

\bibliography{test1}

\end{document}